\documentclass[a4paper,12pt,twoside]{article}

\usepackage[T2A]{fontenc}
\usepackage[utf8]{inputenc}   
\usepackage[russian]{babel}
\usepackage{amsmath,amstext,amsfonts,amsthm,amssymb}
\usepackage{eucal}
\usepackage{graphicx}

\usepackage{fancyhdr,lipsum}
\renewcommand{\subsubsection}[1]{\addtocounter{subsubsection}{1}
{\ \\[3pt]\bf \thesubsubsection. \  #1} }

\swapnumbers

{  \theoremstyle{definition}

}
\newcommand{\Ad}{\operatorname{Ad}}
\newcommand{\ad}{\operatorname{ad}}

\newcommand{\he}{\operatorname{ht}}

\newcommand{\iso}{\overset{\sim}{\longrightarrow}}
\newcommand{\isom}{\overset{\sim}{=}}
\newcommand{\lra}{\longrightarrow}

\newcommand{\dpar}{\partial}

\newcommand{\bea}{\begin{eqnarray*}}
\newcommand{\eea}{\end{eqnarray*}}
\newcommand{\bean}{\begin{eqnarray}}
\newcommand{\eean}{\end{eqnarray}}

\newcommand{\tE}{\tilde E}

\newcommand{\tfh}{\tilde{\mathfrak h}}

\newcommand{\tM}{\tilde M}
\newcommand{\tmu}{\tilde\mu}

\newcommand{\tv}{\tilde v}
\newcommand{\tx}{\tilde x}

\newcommand{\hA}{\hat{A}}

\newcommand{\fg}{\mathfrak g}
\newcommand{\fh}{\mathfrak h}
\newcommand{\fk}{\mathfrak k}

\newcommand{\fgl}{\mathfrak{gl}}

\newcommand{\fsl}{\mathfrak{sl}}

\newcommand{\CD}{\mathcal{D}}

\newcommand{\CL}{\mathcal{L}}

\newcommand{\BC}{\mathbb{C}}

\newcommand{\BR}{\mathbb{R}}
\newcommand{\BZ}{\mathbb{Z}}

\newcommand{\nc}{\newcommand}

\nc{\Id}{\text{Id}}
\nc{\la}{\lambda}

\title{Coxeter element and particle masses}
\author{Laura Brillon\footnote{Institut de Math\'ematiques de Toulouse, 
		Universit\'e Paul Sabatier, 31062 Toulouse, France, Laura.Brillon@math.univ-toulouse.fr} \hspace{1mm} and Vadim Schechtman\footnote{Institut de Math\'ematiques de Toulouse, 
		Universit\'e Paul Sabatier, 31062 Toulouse, France, schechtman@math.ups-tlse.fr}}
\date{22 Novembre 2016}

\begin{document}

\maketitle



\textit{To Joseph Bernstein on his 70th birthday}


\paragraph{Abstract. }
Let $\fg$ be a simple Lie algebra of rank $r$ over $\BC$, $\fh\subset \fg$
a Cartan subalgebra. We construct a family of $r$ commuting
Hermitian operators acting on $\fh$ whose eigenvalues are equal
to the coordinates of the eigenvectors of the Cartan matrix of $\fg$. 

\paragraph{(2010) Mathematics Subject classification.} 17B20, 17B22, 37K10


\section*{Introduction}

This note arose from our attempt to understand a theorem  discovered
by the physicists, [BCDS], [Fr] (a), [FLO], to the effect that the masses of particles (the first $r$ excitations)  
in affine Toda field theories are equal to the coordinates of the Perron -- Frobenius eigenvector of the Cartan matrix $A$. In the text below we review the proof of this elegant result (together with a little generalization), and write down differential equations, similar to 
the Toda field equations, giving rise to particles whose masses 
are absolute values of the coordinates of all other eigenvectors of $A$. 
One observes some interesting regularities in their shape 
related to the geometry of the action of the Coxeter element 
on the Cartan algebra.     

Let $A = (\langle\alpha_i, \alpha_j^\vee\rangle)_{i,j = 1}^r$ be the Cartan matrix of the  root 
system $R\subset \fh^*$ corresponding to a simple finite dimensional complex Lie algebra $\fg$ with a 
fixed Cartan subalgebra $\fh$. 
Let $k_1 < k_2 < \ldots < k_r$ be the exponents of $R$ . Here and below we suppose for simplicity that $R$ is not of type $D_{2n}$, to avoid the case 
when one of the exponents has multiplicity $2$.  
The eigenvalues of $A$ are 
\begin{equation*}
\lambda_i = 2(1 - \cos(k_i\theta)),\ \theta = \pi/h,
\end{equation*}
where $h$ is the Coxeter number of $R$. Let $*:\ \fg\iso\fg$ 
be a Cartan antilinear involution and $H(x,y) = (x, y^*)$ the corresponding Hermitian form, $(x,y)$ being the Killing form. 

The principal result of this paper (see Theorem 5.2) 
is a construction of certain  family of mutually commuting hermitian operators
$M^{(i)}:\ \fh \lra \fh,\ 1\leq i \leq r$, 
such that for each $i$ the eigenvalues 
$\mu^{(i)}_1, \ldots, \mu^{(i)}_r$
of $M^{(i)}$
(in the appropriate order) form an eigenvector of  $A$ with eigenvalue $\lambda_i$.

Actually, the definition of these operators is quite simple 
(we give a sketch here in the introduction, the details the reader will find in the main body of the paper). One starts with a {\it cyclic element} in the sense of Kostant, [K],
\begin{equation*}
e = \sum_{i=0}^r c_ie_i,\ c_i\neq 0,
\end{equation*}
where $e_i\in \fg_{\alpha_i}$, $\alpha_0 := - \theta$,  
$\theta$ being the highest root. Its centralizer $\fh' := Z(e)$ 
is a Cartan subalgebra which is, as Kostant puts it, {\it in apposition} to $\fh$. Let $\fg = \oplus_{i=0}^{h-1} \fg_i$ be the principal gradation, (cf. 1.2). The spaces $\fh^{\prime(i)} := 
\fh'\cap \fg_{k_i}$, $1\leq i\leq r$, are one-dimensional. 

Let $e^{(i)}\in \fh^{\prime(i)}$ be a nonzero vector, for example 
$e^{(1)} = e$. The operators $\ad_{e^{(i)}}\ad_{e^{(i)*}}$ preserve 
$\fh$; let $\tM^{(i)}$ denote its restriction to $\fh$. 
By definition $M^{(i)}$ is a suitable square root of $\tM^{(i)}$. 

The relation to eigenvectors of $A$ is based on a wellknown  relation between $A$ and the Coxeter transformation $c$, due to Coxeter, cf. [Co], (1.5), (1.7),  see \S 3 below for the details\footnote{
Note that  $A$ is (close to) a symmetric matrix, whereas $c$ is an orthogonal matrix, 
the passage from one to another is somewhat similar to the classical  
Cayley transform.}.

The Coxeter element plays a ubiquitous 
role throughout  various domains of  Representation theory, cf. [BGP]. 

An eigenvector $p\in\fh^*$ with the lowest eigenvalue $\lambda_1$, a {\it Perron -- Frobenius vector}, plays a distinguished role. The assertion that the eigenvalues of  $M^{(1)}$ coincide with its coordinates has been  
proven in [Fr] (a), [FLO]; a generalization to $i > 1$ is straightforward. The coordinates of $p$ have a remarkable physical 
meaning, cf. [Cor]. There exist some mysterious formulas expressing them as certain 
products of  values of  Gamma function, cf. [CAS].

In the \emph{Goddard -- Nuyts -- Olive dual} picture the numbers 
$\mu_j^{(i)}$  appear as the charges of static soliton solutions of Toda field 
equations with the purely imaginary coupling constant, corresponding to the Langlands dual Lie algebra $\fg^\vee$, cf. [Fr] (b). 

In the last \S 7 we describe
some  factorization patterns 
in the shape of the Cartan eigenvectors. Namely, among them 
there are exactly $\phi(h)$ vectors \emph{of PF type}  whose 
coordinates are, up to signs, permutations of the coordinates of the PF eigenvector. 
The nonzero components of the other eigenvectors consist of  several  \emph{clusters}, 
each cluster corresponding to a PF eigenvector of a root subsystem $R'\subset R$ with the  Coxeter number $h'$ dividing $h$.

We are grateful to M. Finkelberg for stimulating discussions, and to 
P. Dorey and V. Kac for a useful correspondence. We understand that P. Dorey 
has obtained some results close to our Theorem 5.2.   
Our special thanks go to the referee 
for the useful criticism which allowed us to improve the exposition. 

\section{Principal element and Cartan subalgebras in apposition}

\textbf{1.1. Setup.} Let $\fg$ be a simple finite-dimensional complex  Lie algebra; $( , )$ will denote the Killing form on $\fg$. 
We fix a Cartan subalgebra  
$\fh\subset \fg$; let  $R\subset \fh^*$ be the 
root system of $\fg$ with respect to $\fh$, $\{\alpha_1,\ldots, \alpha_r\}\subset R$ 
a base of simple roots, 
\begin{equation*}
\fg = (\oplus_{\alpha < 0} \fg_\alpha) \oplus \fh \oplus (\oplus_{\alpha > 0} \fg_\alpha),
\end{equation*}
the root decomposition. For $\alpha = \sum_{i=1}^r m_i\alpha_i$ we set 
\begin{equation*}
\he \alpha = \sum_{i=1}^r m_i.
\end{equation*}

Let
\begin{equation*}
\theta = \sum_{i=1}^r n_i\alpha_i ,
\end{equation*}
be the longest root; we set 
\begin{equation*}
\alpha_0 := - \theta, n_0 := 1.
\end{equation*}

The number 
\begin{equation*}
h = \sum_{i=0}^r n_i = 1 + \he \theta ,
\end{equation*}
is the Coxeter number of $\fg$; set $\zeta = \exp(2\pi i/h)$.  

For each $\alpha\in R$ choose a base vector $E_\alpha\in \fg_\alpha$. 

The Killing form 
$( , )$  induces a $W$-invariant scalar product on $\fh$. Using it 
we identify $\fh\iso \fh^*$, and each root may be considered 
as an element of $\fh$. Thus, 
\begin{equation*}
[h, x] = (h, \alpha)x,\ \alpha\in R,\ x\in \fg_\alpha.
\eqno{(1.1.1)}
\end{equation*}

For $x\in \fg$, $Z(x) = \fg^x$ will denote the centralizer of $x$.

Let $G$ denote the adjoint group of $\fg$, e.g., the (Zarisky closure of the) subgroup $G\subset GL(\fg)$ generated by the elements $e^{\ad_x},\ x\in \fg$,   
and 
\begin{equation*}
\exp:\ \fg\lra G 
\end{equation*}
the exponential map. For $g\in G$ and $x\in\fg$, the result 
of the action of $g$ on $x$ will be denoted 
$\Ad_g(x)$. If $\fg$ is realized as a Lie subalgebra of 
a matrix algebra then
\begin{equation*}
\Ad_{\exp(y)}(x) = e^{y}xe^{-y}.
\eqno{(1.1.2)}
\end{equation*}


\textbf{1.2. Principal element and principal gradation.}
Let $\rho^\vee\in \fh$ be defined by
\begin{equation*}
\langle \alpha_i, \rho^\vee\rangle = 1,\ i = 1,\ldots, r.
\end{equation*}

Another definition of $\rho^\vee$:
\begin{equation*}
\rho^\vee = \frac{1}{2}\sum_{\alpha^\vee \in R^\vee_{> 0}}\ \alpha^\vee
\end{equation*}
where $R^\vee\subset \fh$ is the dual root system. It follows that for all $\alpha\in R$ 
\begin{equation*}
\langle \alpha, \rho^\vee \rangle  = \he\alpha.
\end{equation*}

We set 
\begin{equation*}
P = \exp(2\pi i\rho^\vee/h)\in G.
\end{equation*}

For all $\alpha\in R$,
\begin{equation*}
\Ad_P(E_\alpha) = \zeta^{\he\alpha} E_\alpha.
\end{equation*}

Thus, $\Ad P$ defines a $\BZ/h\BZ$-grading on $\fg$, 
\begin{equation*}
\fg = \oplus_{k=0}^{h-1}\ \fg_k,\ \fg_k = \{x\in\fg|\ \Ad_P(x) = 
\zeta^k x\}.
\end{equation*}

We have $\fg_0 = \fh$, and $\fg_1$ admits as a base the set
\begin{equation*}
E_{\alpha_0}, E_{\alpha_1},\ldots, E_{\alpha_r}.
\end{equation*} 

(Note that $\he \alpha_0 = - \he \theta = 1 - h$, 
so that $\Ad_P(E_{\alpha_0}) = \zeta E_{\alpha_0}$.) 


\textbf{1.3. The Cartan subalgebra $\fh'$.} Fix complex numbers $m_i, m'_i$ such that $m_im'_i = n_i$, $i = 0, \ldots, r$, $m_0 = m'_0 = 1$ and define elements 
\begin{equation*}
E = \sum_{i=0}^r m_iE_{\alpha_i}, \tE = \sum_{i=0}^r 
m'_iE_{-\alpha_i}.
\end{equation*}
We have $E\in \fg_1,\ \tE \in \fg_{h-1}$.

{\bf 1.3.1. Lemma.} $[E, \tE] = 0$.

Kostant calls $E, \tE$ {\it cyclic elements}; these are $z_0, \tilde z_0$ in the notation of [K], Thm. 6.7.  

We define, with Kostant, [K], the subspace 
\begin{equation*}
\fh' := Z(E) = Z(\tE) \subset \fg  .
\end{equation*}

It is proven in [K], Thm. 6.7, that $\fh'$ is a Cartan subalgebra of 
$\fg$, called {\it the Cartan subalgebra in apposition to $\fh$ 
with respect to the principal element $P$}. 

The subspace $\fh'\cap\fg_i$ is nonzero iff $i\in \{k_1, k_2,\ldots, 
k_r\}$ where $1 = k_1 < k_2 < \ldots < k_r = h - 1$ are the 
{\it exponents} of $\fg$. We have $k_i + k_{r + 1 - i} = h$. 

Set 
\begin{equation*}
\fh^{\prime(i)} := \fh'\cap\fg_{k_i},\ 1\leq i\leq r;
\end{equation*}
these are the subspaces of dimension $1$.

We denote by 
\begin{equation*}
T' = \exp(\fh')\subset G ,
\end{equation*}
the maximal torus corresponding to $\fh'$.

If $x\in \fh'$, that is, $[x, E] = 0$, then 
\begin{equation*}
0 = [\Ad_P(x), \Ad_P(E)] = \zeta [\Ad_P(x), E],
\end{equation*}
whence $\Ad_P(\fh')\subset \fh'$, e.g., $P\in N_{T'}$. 

Let $R'\subset \fh^{\prime *}$ be the set of roots of $\fg$ 
with respect to $\fh'$, and denote by $\fh^{\prime *}_\BR \subset \fh^{\prime *}$  
the real linear subspace generated by $R'$. 

Recall that the set of (unordered) bases of simple roots in $R'$ is in bijection 
with the set of chambers, the connected components of 
$$
\fh^{\prime *}_\BR\setminus \cup_{\alpha'\in R'} \alpha^{\prime\perp}.
$$
The set of bases is a $W'$-torsor where $W'$ is the Weyl group of $R'$.

A {\it Coxeter element} in $W'$  is an 
element of the form 
$$
c = s_{\alpha'_1}\ldots s_{\alpha'_r}
$$
where $\{\alpha'_1,\ldots,\alpha'_r\}\subset R'$ is some base 
of simple roots. All Coxeter elements are conjugate, cf. [B], Ch. V, 
\S 6, Prop. 1.

\textbf{1.4. Theorem.} (Kostant) The image of $P$ in 
$N_{T'}/T' = W'$ is a Coxeter element.

\begin{proof}[\textbf{Proof.}]
	{\it See [K], Corollary 8.6.}
\end{proof}

As Kostant shows, one can go backwards: 
starting from a Cartan subalgebra $\fh'$ and from an arbitrary 
Coxeter element $c\in W(\fh')$, one can reconstruct $\fh$. We shall use this in \S 5 below.  

\textbf{1.5.} Thus we have 
\begin{equation*}
\fh = \fg^P := \{x\in \fg|\ \Ad_P(x) = x\} ,
\end{equation*}
and
\begin{equation*}
\fh' = \fg^E := \{x\in \fg|\ \Ad_E(x) = [E, x] = 0\}.
\end{equation*}

So, 
\begin{equation*}
E\in \fh'\cap \fg_1,\ \tE\in \fh'\cap \fg_{h-1}.
\end{equation*}
It follows  that $E$ (resp. $\tE$) is an eigenvector of $c$ with eigenvalue $\zeta$ (resp. $\zeta^{-1}$).  
 
\section{Diagonalization of some operators}

\textbf{2.1.} Consider the root decomposition of $\fg$ with respect to $\fh'$: 
\begin{equation*}
\fg = \fh'\oplus_{\alpha'\in R'} \fg_{\alpha'},\ R'\subset \fh^{\prime *}.
\end{equation*}

Recall that $\Ad_P$ leaves $\fh'$ stable and induces the action 
of a Coxeter element $c\in W'$ on $\fh'$, 
\begin{equation*}
W' = N_{T'}/T' \isom N_{\fh'}/\fh'.
\end{equation*}

Recall the the order of $c$ in $W'$ is equal to $h$. 


{\bf 2.2. Proposition.} It is possible to choose nonzero root vectors 
$e_{\alpha'}\in \fg_{\alpha'}$ in such a way that 
$$
\Ad_P(e_{\alpha'}) = e_{c(\alpha')}
$$
for all $\alpha'\in R'$.

\begin{proof}[\textbf{Proof.}]
	{\it See [K], Theorem 8.4. }
\end{proof}

\textbf{2.3.} According to [K], the action of $c$ on $R'$ has $r$ orbits $\Omega_i,\ i = 1,\ldots, r$, each of them containing $h$ elements: 
\begin{equation*}
R' = \coprod_{i=1}^r\ \Omega_i ;
\end{equation*}

here the prime reminds us that the orbits lie in $R'\subset\fh^{\prime *}$.

By the way, it follows that $|R'| = \dim\fg - r = hr$, whence 
\begin{equation*}
\dim \fg = h(r+1).
\end{equation*}

For example, for $\fg = \frak{sl}(n)$, $r = n - 1$, $h = n$ and $\dim \fg = n^2 - 1$. 

For every $1\leq i \leq r$, define with Kostant an element 
\begin{equation*}
a_i = \sum_{\alpha'\in \Omega_i} e_{\alpha'}.
\end{equation*}

It follows from Prop. 2.2 that  
$\Ad_P(a_i) = a_i$, e.g., all $a_i\in \fh = \fg^P$.  

According to [K], Theorem 8.4, $a_1, \ldots, a_r$ forms a base of $\fh$.  

Let us pick an element $\gamma_i\in\Omega_i$, so that 
$\Omega_i = \{c^k(\gamma_i)| k = 0,\ldots, h-1\}$. 

For any  $x\in \fg_m\cap \fh'$, $m\in \BZ/h\BZ$ we have (cf. [Fr] (a))
\begin{equation*}
[x, e_{c^k(\gamma_i)}] = \zeta^{-km}\langle\gamma_i, x\rangle e_{c^k(\gamma_i)}.
\eqno{(2.3.1)}
\end{equation*}
Indeed,
\begin{align*}
[x, e_{c^k(\gamma_i)}] &= [x, \Ad_P^k(e_{\gamma_i})] = 
\Ad_P^k[\Ad_P^{-k}(x), e_{\gamma_i}] \\
&= \zeta^{-km}\Ad_P^k[x, e_{\gamma_i}] = \zeta^{-km}\Ad_P^k(\langle\gamma_i, x\rangle e_{\gamma_i})\\
&= \zeta^{-km}\langle\gamma_i, x\rangle e_{c^k(\gamma_i)}. 
\end{align*}

It follows that for any $y\in \fg_{-m}\cap\fh'$,
\begin{equation*}
[y,[x, e_{c^k(\gamma_i)}]] = \langle\gamma_i, y\rangle \langle\gamma_i, x\rangle e_{c^k(\gamma_i)}.
\end{equation*}

Summing up by $k$, we get the following theorem.

{\bf 2.4. Theorem.}  For all $m\in \BZ/h\BZ$, $x\in \fg_m\cap \fh',\ y\in \fg_{-m}\cap\fh'$, $1\leq i\leq r$, 
$$ 
[y, [x, a_i]] = \langle\gamma_i, y\rangle \langle\gamma_i, x\rangle a_i.
$$
In other words, $\{a_1,\ldots, a_r\}$  is  a base of $\fh$ which diagonalizes the operator $\ad_y\ad_x$.

\section{Coxeter element, Cartan matrix, and their   eigenvectors  }

\textbf{3.1.} Let $R\subset V$ be a reduced irreducible root system  in a real vector space $V$ of dimension $r$ (in particular $R$ generates $V$), $W$ the Weyl group of $R$, $( , )$ a $W$-invariant scalar product on $V$. We identify $V$ with $V^*$ using $( , )$ so that $\alpha^\vee = 2\alpha/(\alpha, \alpha)$, cf. [B], Chapitre V, \S 1, Lemme 2.   

Let $\{\alpha_1, \ldots, \alpha_r\} \subset R$ be a base of simple roots.

Choose a  black/white {\it colouring} of the set $I$ of vertices of the Dynkin graph of $R$ (which is a tree) 
such that neighbouring vertices have different colours. Identify $I$ with $\{1,\ldots, r\}$ in such a way that the vertices $\{1,\ldots, p\}$ are black, and the vertices $\{p+1,\ldots, r\}$ are white. 

We denote $s_i := s_{\alpha_i}$. Consider a Coxeter element
\begin{equation*}
c = c_bc_w,\ c_b = \prod_{i=1}^p s_i,\ c_w = \prod_{i=p+1}^r s_i,
\end{equation*} 
the order inside the products defining $c_b$ and $c_w$ does not matter since reflections $s_i, s_j$ commute once $i$ and $j$ have the same colour. Evidently 
\begin{equation*}
c_w^2 = c_b^2 = 1,
\end{equation*}
whence
\begin{equation*}
(c_b + c_w)^2 = c + c^{-1} + 2.
\end{equation*}

Let $A = (n_{ij}) = (\langle \alpha_i, \alpha_j^\vee\rangle )_{i,j = 1}^r $ be the Cartan matrix of $R$. We denote by $\hat A: V\lra V$ an operator defined by 
\begin{equation*}
\hat A(\alpha_i) = \sum_{j=1}^r n_{ij}\alpha_j.
\end{equation*}

{\bf 3.2. Lemma.}  We have 
\begin{equation*}
c_b + c_w = 2I - \hat A. 
\end{equation*}

\begin{proof}[\textbf{Proof.}]
{\it The matrix $A$ has a block form
\begin{equation*}
A = \left(\begin{matrix} 2I_p & X\\ Y & 2I_{r-p}\end{matrix}\right).
\end{equation*}
	
On the other hand, calculating the action of the (commuting) simple reflections, one finds the matrices of the operators $c_b$ and $c_w$ in the base $\alpha_1, \ldots, \alpha_r$ to be 
\begin{equation*}
c_b = \left(\begin{matrix} -I & -X\\ 0 & I\end{matrix}\right),
c_w = \left(\begin{matrix} I & 0\\ -Y & -I\end{matrix}\right), 
\end{equation*}
whence 
\begin{equation*}
c_b + c_w = \left(\begin{matrix} 0 & -X\\ -Y & 0\end{matrix}\right) = 2I - \hat A.
\end{equation*} }
\end{proof}

{\bf 3.3. Lemma.} All the eigenvalues of $\hA$ have the form $2(1 - \cos k_i\theta_1)$, $i \in \{1,\ldots, r\}$ where $k_i$ are the exponents of $\fg$ and $\theta_1 = \pi/h$.

\begin{proof}[\textbf{Proof.}]
{\it We shall use the identity
\begin{equation*}
(2I - \hA)^2 = c + c^{-1} + 2.
\eqno{(3.3.1)}
\end{equation*}

The eigenvalues of $c$ are $e^{2k_j\pi i/h}, 1\leq j\leq r$, cf. [Co]. It follows from (3.3.1) that if $e^{2i\theta}$ is an eigenvalue of $c$, then $4\cos^2\theta$ is an eigenvalue of $(2I - \hA)^2$, so $2(1 \pm \cos\theta)$ is an eigenvalue of $\hA$. 
		
Note that 
\begin{equation*}
2(1 + \cos\theta) = 2(1 - \cos(\pi - \theta)),
\end{equation*}
and  $k_{r - i} = h - k_i$, which implies the assertion of the lemma.}
\end{proof}

\textbf{3.4.} For a vector 
\begin{equation*}
x = \sum_{i=1}^r x_i\alpha_i\in V,
\end{equation*}
we have 
\begin{equation*}
\hat A x = \sum_i x_i( \sum_j n_{ij})\alpha_j = 
\sum_j (\sum_ix_in_{ij})\alpha_j  .
\end{equation*}
Thus, $\hat A x = \lambda x$ is equivalent to 
\begin{equation*}
\sum_i\ x_in_{ij} = \lambda x_j,\ j = 1,\ldots, r.
\eqno{(3.4.1)}
\end{equation*}

Define a {\it colour function}\ 
$\epsilon: \{1,\ldots, r\} \lra\{\pm 1\}$ by 
\begin{equation*}
\epsilon(i) = \left\{\begin{matrix} 1 & \text{if}\ i\leq p \\
-1 & \text{if}\ i > p .\\
\end{matrix}\right.
\end{equation*}

{\bf 3.5. Duality Lemma.}  Let $x$ satisfy (3.4.1)  
with $\lambda = 2(1 - \cos\theta)$, that is,
\begin{equation*}
\sum_{i=1}^r\ x_in_{ij} = 2(1 - \cos\theta)  x_j,\ j = 1,\ldots, r.
\eqno{(3.5.1)}
\end{equation*}

Then 
\begin{equation*}
\sum_{i=1}^r\ \epsilon(i)x_in_{ij} = 2(1 + \cos\theta) \epsilon(j) x_j,\ j = 1,\ldots, r.
\eqno{(3.5.2)}
\end{equation*}

\begin{proof}[\textbf{Proof.}]
{\it Recall that the matrix $A$ has a block form $\begin{pmatrix} 2 I_p & X \\ Y & 2 I_{r-p} \end{pmatrix}$ and write $x_b = ( x_1 , .. , x_p)$ and $x_w = (x_{p+1} , .. , x_r ) $. The identity $xA = \lambda x$ gives :
	
\begin{equation*}
\left\{
\begin{array}{l}
2 x_b + x_w Y = 2 ( 1 - \cos \theta ) x_b \\
x_b X + 2 x_w = 2 ( 1 - \cos \theta ) x_w .
\end{array}
\right.
\end{equation*} 
Then
\begin{equation*}
\left\{
\begin{array}{l}
2 x_b - x_w Y = 2 ( 1 + \cos \theta ) x_b \\
x_b X - 2 x_w = 2 ( 1 + \cos \theta ) ( - x_w ).
\end{array}
\right.
\end{equation*}
This means that $\tilde{x} =  ( \epsilon(1)x_1 , ... , \epsilon(r)x_r )$ satisfies (3.4.1) with $\lambda = 2 (1+ \cos \theta )$. }	
\end{proof}

\subsection{ }
Now set 
\begin{equation*}
x_b = \sum_{i=1}^p x_i\alpha_i,\ 
x_w = \sum_{i=p+1}^{r} x_i\alpha_i.
\end{equation*}

{\bf Lemma.} (a) $c_w(x_w) = - x_w,\ c_b(x_b) = - x_b$.
(b) $c_w(x_b) = x_w + 2\cos\theta x_b, c_b(x_w) = x_b + 2\cos\theta x_w$.


{\bf 3.7. Corollary.} Define $y = e^{-i\theta/2}x_w + e^{i\theta/2}x_b$. Then
\begin{equation*}
c(y) = e^{2i\theta} y.
\end{equation*}


{\bf 3.8. Lemma.} For all $j = 1, \ldots, r$, 
\begin{equation*}
 (y, \alpha_j)  = i\epsilon(j)e^{-i\epsilon(j)\theta/2}\sin \theta\cdot (\alpha_j,\alpha_j) x_j.
\end{equation*}


\begin{proof}[\textbf{Proof.}]
{\it Recall that $(\alpha_k , \alpha_j) = \frac{1}{2} (\alpha_j , \alpha_j)n_{kj}, k,j = 1 , \ldots , r$. We have
\begin{equation*}
( y , \alpha_j ) = \frac{ (\alpha_j , \alpha_j )}{2} \left( e^{i\theta/2} \sum_{k=1}^{p} x_k n_{kj} + e^{-i \theta/2} \sum_{k=p+1}^{r} x_k n_{kj} \right) .
\end{equation*} 

Since
\begin{equation*}
\sum_{k=1}^{p} x_k n_{kj} = \frac{1}{2} \left( \sum_{i=1}^{r} x_i n_{ij} + \sum_{i=1}^{r} \epsilon(i)x_i n_{ij} \right),
\end{equation*}
and
\begin{equation*}
\sum_{k=p+1}^{r} x_k n_{kj} = \frac{1}{2} \left( \sum_{i=1}^{r} x_i n_{ij} - \sum_{i=1}^{r} \epsilon(i)x_i n_{ij} \right),
\end{equation*}
the application of Lemma 3.5 gives the result.}	
\end{proof}

{\bf 3.9. Lemma.} The elements $\epsilon(i)\alpha_i,\ 1\leq i\leq r$, belong to $r$ different orbits of the action of $c$ on $R$.

\begin{proof}[\textbf{Proof.}]
{\it Kostant proves in [K], Thm 8.1 and Thm 8.4, that exactly $r$ negative roots, say $\{ \beta_1 , ... , \beta_r \}$, become positive under the action of $c$ on $R$ and they belong to $r$ different orbits of this action.
	
In the proof of Lemma 3.2, we have seen that 
\begin{equation*}
 c = c_b c_w = \begin{pmatrix}
 -I + XY & X \\
 -Y & -I \end{pmatrix} \text{ and } c^{-1} = c_w c_b = \begin{pmatrix}
 -I & -X \\
 Y & YX - I \end{pmatrix} ,
\end{equation*} 
where $X$ and $Y$ are matrices with nonpositive entries, such that $A = \begin{pmatrix}
2 I_p & X \\
Y & 2 I_{r-p} \end{pmatrix} $. 

For $1 \leq i \leq p$, $c^{-1}(\alpha_i)$ is a negative root, whence $\alpha_i =c(\beta_k)$ for some $k$ in $\{ 1 , .. , r\}$. For $p+1 \leq i \leq r$, $c(-\alpha_i)$ is a positive root, whence $- \alpha_i = \beta_j$ for some $j$ in $\{1,..,r\}$. 

Thus, to each root $c(i)\alpha_i$, we have associated a root $\beta_j$ in the same orbit, and this is a one-to-one correspondence.
}	
\end{proof}

\section{Cartan involution and  Hermitian form}

\textbf{4.1.} Recall the setup 1.1. Let us choose, with F.Bruhat [Br] 
and Kostant [K], p. 1003,  
a {\it Weyl basis} $\{ e_\alpha\in\fg_\alpha\}$. By definition, this means that all $e_\alpha\neq 0$, $(e_\alpha, e_{-\alpha}) = 1$, and 
if we denote
\begin{equation*}
[e_\alpha, e_\beta] = n_{\alpha\beta} e_{\alpha + \beta},
\end{equation*}
then $n_{\alpha\beta} = n_{-\alpha,-\beta}$. Here we have chosen the Kostant's normalization of the Weyl basis. 
We set $h_i := [e_{\alpha_i}, e_{-\alpha_i}],\ 1\leq i\leq r$.   

Let $\fk\subset \fg$ be the real Lie subalgebra with the base
\begin{equation*}
e_\alpha - e_{-\alpha},\ i(e_\alpha + e_{-\alpha}),\ ih_j , \ \alpha\in R_+,\ 1\leq j\leq r .
\end{equation*}

It is a {\it compact form} of $\fg$, which means by definition that
\begin{equation*}
\fg = \fk\oplus i\fk,
\end{equation*}
and the restriction of the Killing form to $\fk$ is negative definite. 


Define, following Kostant, an
involution $(\cdot)^*: \fg\iso \fg$ by 
$$
(x + iy)^* = x - iy,\ x, y\in i\fk.
$$
Then 
$(\lambda x)^* = \bar \lambda x, \lambda\in\BC$, and 
$$
[x, y]^* = [y^*, x^*].
$$
The sesquilinear form on $\fg$
$$
H(x, y) = (x, y^*)
$$
is Hermitian positive definite, cf. [Br], (21). 

With respect to this form
$$
(\ad_x)^* = \ad_{x^*}
$$
One computes that
$$
e_{\alpha}^* = e_{-\alpha},\ \alpha\in R.
\eqno{(4.1.1)}
$$

\textbf{4.2.} Now let us return to the setup of Sections 1 and 2. 
However now we will work with specifically chosen Weyl vectors $e_\alpha$, instead of arbitrary root vectors 
$E_\alpha$. 

Fix nonzero complex numbers $m_i$ such that 
$m_i\bar m_i = n_i$, $1\leq i\leq r$. 
Let 
$$
e = \sum_{i=1}^r m_i e_{\alpha_i} + e_{-\theta}
$$
be the cyclic element. 

By (4.1.1),  
$$
e^* = \sum_{i=1}^r \bar m_i e_{-\alpha_i} + e_{\theta}.
$$
Recall that $[e,e^*] = 0$.

Let $\fh' = Z(e) = Z(e^*)$ be the corresponding Cartan subalgebra 
in apposition to $\fh$, as in 1.3. 

For all $m\in \BZ/h\BZ$
$$
*: \fg_m\cap \fh' \iso \fg_{-m}\cap \fh' 
$$

\textbf{4.3.} Let us apply Theorem 2.4 to 
$y = x^*$. 
With $a_i\in \fh$ and $\gamma_i\in R'\subset \fh^{\prime *}$ as in 2.3,  we obtain
$$
\ad_x\ad_{x^*}(a_i) = \gamma_i(x)\gamma_i(x^*)a_i,\ 1\leq i\leq r.
$$

\textbf{4.4. }{\bf Lemma.}  $\gamma_i(x^*) = \overline{\gamma_i(x)}$.

\begin{proof}[\textbf{Proof.}]
{\it Consider the equality (2.3.1): 
	$$
	\ad_x(z) = \zeta^{-km}\gamma_i(x)z
	$$
	where we set for brevity $z = e_{c^k(\gamma_i)}$. It follows that 
	$$
	H(\ad_x(z), z) = \zeta^{-km}\gamma_i(x)H(z,z).
	$$
	Similarly, 
	$$
	\ad_{x^*}(z) = \zeta^{km}\gamma_i(x^*)z,
	$$
	whence 
	$$
	H(\ad_{x^*}(z), z) = \zeta^{km}\gamma_i(x^*)H(z,z).
	$$
	The assertion follows now from the adjointness of the operators 
	$\ad_x$ and $\ad_{x^*}$, since $H(z,z) \neq 0$. }	
\end{proof}

\section{Main theorem}

\textbf{5.1.} Let us start with a Cartan subalgebra $\fh'\subset \fg$, 
whence the root system $R' \subset \fh^{\prime *}$. Choose a base  
of simple roots $\{\alpha'_i\}\subset R'$ and a bicolouring 
of the Dynkin graph as in \S 3. Thus, 
$\alpha_i'$ with $1\leq i \leq p$ (resp. with $p+1\leq i\leq r$) 
will denote the \emph{black} (resp. \emph{white}) simple roots. 
Let 
$$
c' = c'_bc'_w,\ c'_b = \prod_{i=1}^p s'_i,\ c'_w = \prod_{i=p+1}^r s'_i,
$$
where $s'_i := s_{\alpha'_i}$ be the corresponding Coxeter element. 

Let $G$ be the adjoint group of $\fg$, $T'\subset G$ the 
maximal torus with $\text{Lie}(T') = \fh'$, so that the Weyl group 
$W'\subset \text{Aut}(R')$ will be identified with $N_G(T')/T'$. 
Let $P'\in N_G(T')$ be an element that projects to $c'$. 
Set
$$
\fh = \fg^{P'}.
$$
Then $\fh$ is a Cartan subalgebra, and $\fh'$ is in apposition to 
$\fh$ with respect to $P'$, cf [K], Theorem 8.6\footnote{The couple of Cartan subalgebras $(\fh, \fh')$ from [K], \S 6 becomes $(\tilde\fh, \fh)$ in 
{\it op. cit.}, 8.6. Two occurences of $\fh$ in [K], p. 1023, 2nd line, should be replaced by $\tfh$.}. 

Consider the principal gradation generated by $P'$, 
$\fg = \oplus \fg_i$, where $\fg_i$ is the $\zeta^i$-eigenspace of $\Ad_{P'}$, 
as in \S 1, 
and one-dimensional spaces $\fh^{\prime(i)} := \fh'\cap \fg_{k_i}$, $1\leq i\leq r$. 

We can choose an involution $*$ as in \S 4 in such a way that it leaves $\fh'$ invariant,  for every $x\in \fh'$ the operator $\ad_{x^*}$ is Hermitian 
conjugate to $\ad_{x}$, and $(\fh^{\prime(i)})^* = \fh^{\prime(r-i)}$. 

Indeed, this is true for  the gradation induced 
by the principal element $P = P_0$ as defined in 1.2, and the involution  (let us denote it $*_0$) defined as in \S 4 starting from $\fh$. Afterwards one can use the conjugacy Theorem 7.3 from [K] 
to define the desired involution for the principal gradation induced by $P'$.


{\bf 5.2. Theorem.} Let $i$ be an integer, $1\leq i\leq r$. Let $e^{(i)}$ be a nonzero vector in $\fh^{\prime(i)}$,  whence  $e^{(i)*}\in \fh^{\prime(r - i)}$. Consider 
a selfadjoint nonnegative  operator 
$$
\tM^{(i)} := \ad_{e^{(i)}}\ad_{e^{(i)*}}:\ \fh\lra \fh.
$$
Let
$\tmu_1^{(i)},\ldots,\tmu_r^{(i)}$ 
denote its eigenvalues.

There exists an (essentially unique) operator $M^{(i)}\in \fgl(\fh)$ whose square is equal to $\tM^{(i)}$ such that the column vector of its eigenvalues in the approriate numbering 
$$
\mu^{(i)} := 
(\mu_1^{(i)},\ldots,\mu_r^{(i)})^t
\eqno{(5.2.1)}
$$
is an  eigenvector of  the Cartan matrix $A$ with eigenvalue 
$$
\lambda_i := 2(1 - \cos(2k_i\pi/h)).
$$ 
 In particular, for $i = 1$ there exists an eigenvector of $A^t$ 
with all coordinates positive (a {\bf Perron -- Frobenius vector}),  
and we may take as $M^{(1)}$ the positive square root of $\tM^{(1)}$.

The operators $M^{(1)}, \ldots, M^{(r)}$ commute with each other.

\begin{proof}[\textbf{Proof.}]
{\it All the necessary tools have been already prepared.  Let 
	$x = (x_1,\ldots, x_r)^t$ be an eigenvector of $A^t$ (sic!) 
	with eigenvalue $\lambda = \lambda_i$. Starting from it, define an eigenvector $y$ of the Coxeter element $c'$, with eigenvalue 
	$e^{\sqrt{-1}\theta}$, $\theta = k_i\pi/h$, cf. Corollary 3.7. By Lemma 3.8, we have
	$$
	(y, \alpha'_j)  = \sqrt{-1}\epsilon(j)e^{-\sqrt{-1}\epsilon(j)\theta/2}\sin \theta\cdot (\alpha'_j,\alpha'_j) x_j, \ 1\leq j\leq r.
	\eqno{(5.2.2)} 
	$$
	
On the other hand, we know from \S 4 the eigenvalues of   
the  operator $\tM = \tM^{(i)}$ : as follows from 4.3 and 4.4, they are 
$|\gamma_j(e^{(i)})|^2, 1\leq j\leq r$
where $\gamma_j\in \fh^{\prime *},\ 1\leq j\leq r$, are arbitrary representatives of different orbits of the action 
of $c$ on $R'$, 

Let us identify $\fh'$ with $\fh^{\prime *}$ 
using the scalar product $( , )$, so that we can consider 
$\gamma_j$ as vectors belonging to $\fh'$, and $|\gamma_j(e^{(i)})|^2 = 
|(e^{(i)},\gamma_j)|^2$.

Recall that $e^{(i)}$ is an eigenvector of 
$c'$ in $\fh'$ : $c'(e^{(i)}) = \lambda e^{(i)}$ ($c'$ acts as $\Ad_{P'}$ on $\fh'$), whence
$$
e^{(i)} = \mu y
$$
for some $\mu\in \BC^*$.

Let us rewrite (5.2.2) in the form
$$
(y, \epsilon(j)\alpha'_j)  = \sqrt{-1}e^{-\sqrt{-1}\epsilon(j)\theta/2}\sin \theta\cdot\tx_j, \ 1\leq j\leq r,
\eqno{(5.2.2)} 
$$
where $\tx_j = (\alpha'_j,\alpha'_j) x_j$. 

Note that 
$\tx = (\tx_1,\ldots,\tx_r)^t$ is a  $\lambda$-eigenvector of $A$. 

Due to Lemma 3.9, the vectors $ \epsilon(j)\alpha'_j, 1\leq j\leq r$, are  representatives 
of $r$ orbits of $c'$-action on $R'$, so we can set 
$$
\gamma_j := \epsilon(j)\alpha'_j.
$$
It follows that the eigenvalues of $\tM = \tM^{(i)}$ 
are 
$$
|(e^{(i)}, \gamma_j)|^2 =  |\mu|^2\sin^2\theta \cdot \tx_j^2,\ 1\leq j\leq r
$$
(note that $\tx_j$ are real numbers, not necessarily positive). 

Thus, the sequence of eigenvalues of $\tM^{(i)}$ is 
$$
(|\mu|^2\sin^2\theta \tx_1^2,\ldots, |\mu|^2\sin^2\theta\tx_r^2).
$$
On the other hand, a $\lambda_i$-eigenvector of $A$ is
$$
(\tx_1,\ldots,\tx_r)^t
$$
Moreover, the operators $\tM^{(1)},\ldots, \tM^{(r)}$ mutually commute since 
the elements $e^{(1)},\ldots, e^{(r)}\in \fh'$ mutually commute. 
	
Now, as an operator 
$M = M^{(i)}$, we take (the unique) one of the $2^r$ square roots of $\tM$ 
whose $j$-th eigenvalue, if nonzero, has the same sign as that of $x_j$.
The set of eigenvalues of $M^{(i)}$ will be 
$$
(|\mu \sin\theta| \tx_1,\ldots, |\mu\sin\theta|\tx_r),
$$
and this vector is a  $\lambda_i$-eigenvector of $A$.	} 
\end{proof}


\textbf{5.3.} We can start with an arbitrary pair of Cartan subalgebras 
$\fh, \fh'$ where $\fh'$ is in apposition to $\fh$ with respect to a principal 
element $P$.  Defining operators $\tM^{(i)}$ as in 5.2, we arrive at the same 
conclusions for their spectra as in 5.2, due again to the Kostant's 
conjugacy theorem, [K], Theorem 7.3. 

\section{Affine Toda field equations }

\textbf{6.1. Affine Toda field theories.} Consider a classical field theory whose 
fields  are smooth functions
$\phi: X\lra \fh$
where $X = \BR^2$ \emph{space - time}, with coordinates $x_1, x_2$.  

The Lagrangian density of the theory depends on an element $e\in\fh'$ 
where  $\fh'$ is a Cartan algebra in apposition to 
$\fh$,  and is given by  
$$
\CL_e(\phi) = \frac{1}{2}\sum_{a = 1}^2 (\dpar_a\phi, \dpar_a\phi) 
- m^2 (\Ad_{\exp(\phi)}(e), e^*).
$$
Here $\dpar_a := \dpar/\dpar x_a$. 

The Euler --Lagrange equations of motion are 
$$
\CD_e(\phi) := 
\Delta\phi + m^2[\Ad_{\exp(\phi)}(e), e^*] = 0,
\eqno{(6.1.1)}
$$
where $ \Delta\phi = \sum_{a=1}^2 \dpar^2_a\phi$. 
It is a system of $r$ nonlinear differential equations 
of the second order. To write them down explicitly 
one uses the formula (1.1.2). 

The usual ATFT correponds to the choice of $e\in \fg^{(1)}$ as in 4.2, cf. [Fr].


The linear approximation to the nonlinear equation (6.1.1) is 
a Klein -- Gordon equation  
$$
\Delta_e\phi :=  \Delta\phi + m^2\ad_e\ad_{e^*}(\phi) = 0 
\eqno{(6.1.2)}
$$
It admits $r$ \emph{normal mode} solutions
$$
\phi_j(x_1, x_2) = e^{i(k_jx_1 + \omega_j x_2)}y_j,\ k_j^2 + \omega_j^2 = m^2\mu_j^2,
$$
$\ 1\leq j\leq r$, where $\mu_j^2$ are the eigenvalues of the 
square mass operator 
$$
M_e^2 := \ad_e\ad_{e^*}:\ \fh\lra \fh
$$
and $y_j$ are the corresponding eigenvectors, cf. [H] (1.4), (1.5). 

In other words, (6.1.2) decouples into $r$ equations describing 
scalar particles of masses $\mu_j$, which explains the name \emph{masses}  for them. 


Due to commutativity of $\fh'$, for 
all $e, e'\in \fh'$, 
$$
[ \Delta_e, \Delta_{e'}] = 0.
$$

\section{Factorization patterns in Cartan eigenvectors}

\textbf{7.1.} We recall the 
eigenvectors $\mu^{(i)}$ from Theorem 5.2; they are in bijection 
with the exponents $k_i$, $1\leq i\leq r$. In particular, 
$\mu^{(1)}$ is a Perron -- Frobenius eigenvector. 

The exponents come in pairs $k_i, k_{r-i} = h - k_i$. According 
to Lemma 3.5, the eigenvector $\mu^{(r-i)}$ is obtained from 
$\mu^{(i)}$  
by multiplying the coordinates by the sequence $(\epsilon(1),\ldots,\epsilon(r))$, 
with $\epsilon(j) = \pm 1$. 

Below we will use the following notation. For a	vector $v = (x_1, \ldots, x_r)$ 
we denote $\tv = (|x_1|, \ldots, |x_r|)$. For $\sigma\in S_r$ 
(the symmetric group), $v_\sigma  := (x_{\sigma(1)},\ldots, x_{\sigma(r)})$.  

The notation $\gcd(a,b)$ will mean the greatest common 
divisor of $a$ and $b$. 

Consider first the case $\fg = \fsl(m)$. A Perron--Frobenius vector for the 
Lie algebra $ \fsl(m)$ has the form
$$
\mu^{(1)} = v_{PF}(m) := (\sin(\pi/m),\ldots,\sin((m-1)\pi/m) 
$$ 

Let  $\fg = \fsl(n)$, and let us describe the other eigenvectors $\mu^{(i)}$, 
$1\leq i\leq r = n-1$. 

Let $p(i) = \gcd(i,n)$, $n = p(i)q(i)$. 

Consider first the case $p(i) = 1$. In this case it is not difficult to see that all the components 
of $\mu^{(i)}$ are nonzero, and form, up to a sign, a permutation 
of the components of $\mu^{(1)}$.

(The permutations involved will be described in 7.2 below. )

For an arbitrary $p(i)$, 
among the components of $\mu^{(i)}$ 
there are exactly $p(i) - 1$ zeros, and
the remaning $p(i)(q(i) - 1)$ components may be decomposed into $p(i)$ groups, the numbers inside each group 
forming, up to a sign, a Perron -- Frobenius 
eigenvector for $\fsl(q(i))$.

{\bf 7.1.1. Example.} For the Cartan matrix of $\fsl(12)$, 
the eigenvectors are $\mu^{(i)}, 1\leq i\leq 11$, with 
$\tmu^{(i)} = \tmu^{(12 - i)}$. Then we have
\newline  
$\tmu^{(2)} = \left( v_{PF}(6) , 0 , v_{PF}(6) \right)$, 
$\tmu^{(3)} = \left( v_{PF}(4) , 0 , v_{PF}(4) , 0 , v_{PF}(4) \right)$
\newline
$ \tmu^{(4)} = \left( v_{PF}(3) , 0 , v_{PF}(3) , 0 , v_{PF}(3) , 0 , v_{PF}(3) \right) $
$\tmu^{(5)} =  v_{PF}(12)_\sigma$, with  $\sigma = ( 1 \text{ } 5)(7 \text{ } 11 )$,
\newline
$ \tmu^{(6)} = \left( v_{PF}(2) , 0 ,v_{PF}(2) , 0 ,v_{PF}(2) , 0 ,v_{PF}(2) , 0 ,v_{PF}(2) , 0 ,v_{PF}(2) \right) $. 
\flushright
 $\square$
\flushleft

For an arbitrary $\fg$ we have a similar pattern. 
  Let $R$ be a finite 
reduced irreducible root system of rank $r$ and the Coxeter number $h$, $1\leq i\leq r$, $k_i$ the corresponding exponent.

Note that according to [B], Chapitre VI, \S 1, Proposition 30, 
all numbers $1\leq k\leq h-1$ prime to $h$ are among the exponents. 

Let $p(i) = \gcd(k_i, h)$. 


{\bf 7.2. Proposition.} Suppose that $R$ is simply laced. (a) The eigenvector $\mu^{(i)}$ has all components 
different from $0$ iff $p(i) = 1$, and if this is the case, we have 
$$
\tilde\mu^{(i)} = \tilde\mu^{(1)}_{\sigma_i}
$$
for some $\sigma_i\in S_r$. 

There are $\phi(h)/2$ such permutations $\sigma_i$, and they form a group isomorphic to $U(\BZ/h\BZ)/\{1, -1\}$.

(b) If $p(i)$ is  arbitrary, then one can associate to such $i$ 
a root subsystem 
$R_{p(i)}\subset R$ whose Coxeter number is $q(i) = h/p(i)$ in such a way that  the nonzero components of $\mu^{(i)}$ 
are decomposed into $p(i)$ groups, each group  being, 
up to signs, a permutation of the coordinates of a PF vector  
for  $R_{p(i)}$.

\vspace{3mm}

These facts may be verified  case-by-case, using the explicit formulas 
for the vectors $\mu^{(i)}$ given in [Do], Table 2 on p. 659.

However, it would be desirable to have a uniform proof of this. 

We believe that the same holds true for non-simply laced $R$ as well.

\textbf{7.3. Example. } For the root system of type $E_8$, we have $h = 30$, 
the exponents are 
\newline $1, 7, 11, 13, 17, 19, 23, 29$; they include exatly  all prime 
numbers $\leq 30$ not dividing $30$ (and $1$).  
Let us denote 
the corresponding eigenvectors $v_1, \ldots, v_{29}$, so that 
$v_k$ has eigenvalue $2(1 - \cos(k\pi/30))$. The first one 
$v_1 = v_{PF}$ is a Perron -- Frobenius vector. It is equal to 

$v_{PF} = $
$$
(1, \frac{1}{\mu} (\mu^2 - 1), 
\mu, \mu^2 - 1, 
\frac{1}{\mu} ( \mu^4 - 3 \mu^2 + 1 ),
\mu^4 - 4 \mu^2 + 2, 
\frac{1}{\mu} ( \mu^6 - 5 \mu^4 + 5 \mu^2 - 1 ), 
\mu^6 - 6 \mu^4 + 9 \mu^2 - 3)
$$
where $\mu = 2\cos(\pi/30)$. 

Then we have 
$\tv_1 = \tv_{29},\ \tv_{7} = \tv_{23}, \tv_{11} = \tv_{19}, 
\tv_{13} = \tv_{17}$
and
$$
\tv_7 = (v_{1})_\sigma, \tv_{11} = (v_1)_{ \sigma^2}, \tv_{13} = 
(v_1)_{ \sigma^3} 
$$
with $\sigma = (1742)(3658)\in S_8$. The cyclic subgroup 
$G = \{1, \sigma, \sigma^2, \sigma^3\}\subset S_8$ is isomorphic 
to $U(\BZ/30\BZ)/\{1, - 1\}$. 
\flushright
$\square$    
\flushleft

\newpage 
    
\section*{References} 

[B] N. Bourbaki, \textit{Groupes et alg\`ebres de Lie.} Chapitres 4, 5 et 6. Paris, 
Hermann, 1968. 

[BGP] I.N. Bernstein, I.M. Gelfand, V.A. Ponomarev, \textit{Coxeter functors and a theorem of Gabriel} ; И.Н.Бернштейн, И.М.Гельфанд, В.А.Пономарев, 
Функторы Кокстера и теорема Габриэля, {\it УМН} {\bf 28} (1973), 19 -- 33.

[BCDS] H.W. Braden, E. Corrigan, P.E. Dorey, R. Sasaki, \textit{Affine Toda field 
theory and exact $S$-matrices},  Nucl. Phys. {\bf B338} (1990), 
689 -- 746.

[Br] F. Bruhat, \textit{Formes r\'eelles des alg\`ebres semi-simples},
 S\'eminaire Sophus Lie, {\bf 1} (1954 - 1955), Expos\'es n$^o$ 11 -- 12. 

[Ca] B. Casselman, \textit{Essays on Coxeter groups.} Coxeter elements in finite 
Coxeter groups, https://www.math.ubc.ca/~cass/research/pdf/Element.pdf.

[CAS] V. Cohen-Aptel, V. Schechtman, \textit{Vecteurs de Perron--Frobenius 
et produits Gamma}, Comptes Rendus Math\'ematique {\bf 350} 
(2012), 1003 -- 1006.


[Cor] E. Corrigan, \textit{Recent developments on affine Toda field theory}, 
in: G. Semenov et al. (eds),  Particles and fields (1999), 
1 -- 34. 

[Co] H.S.M. Coxeter, \textit{The product of the generators of a finite group 
generated by reflections},  Duke Math. J. {\bf 18} (1951), 
765 -- 782.

[D] P. Dorey, \textit{Root systems and purely elastic $S$-matrices}, Nucl. Phys. {\bf B358} (1991), 654 -- 676.

[Fr] M.D. Freeman, (a) \textit{On the mass spectrum of affine Toda field theory}, Phys. Let. B {\bf 261} (1991), 57 -- 61; (b) \textit{Conserved charges and soliton solutions in affine 
Toda theory}, Nucl. Phys. {\bf B433} (1995), 657 -- 670.  

[FLO] A. Fring, H.C. Liao, D.I. Olive, \textit{The mass spectrum and coupling in affine Toda theories}, Phys. Let. B {\bf 266} (1991), 82 -- 86.

[H] T. Hollowood, \textit{Solitons in affine Toda field theories}, 
Nucl. Phys. {\bf B384} (1992), 523 -- 540.

[K] B. Kostant, \textit{The principal three-dimensional subgroup 
and the Betti numbers of a complex simple Lie group}, Amer. J. Math. 
{\bf LXXXI} (1959),  973 -- 1032. 


\end{document}